\documentclass{article}

\usepackage{amsmath}
\usepackage{amssymb}
\usepackage{theorem}
\usepackage{xypic}
\usepackage[all,v2]{xy}
\xyoption{2cell}
\UseAllTwocells
\usepackage{enumitem}
\usepackage{rotating}
\usepackage{mathrsfs}
\usepackage{stmaryrd}
\usepackage{mathdots}
\usepackage{hyperref}
\usepackage{caption}
\newcommand{\widebar}{\ol}

\newenvironment{items}
{\begin{enumerate}[topsep=3pt, itemsep=3pt, parsep=0pt, label=(\roman*)]}
{\end{enumerate}}

\renewcommand{\tilde}{\widetilde}

\newcommand{\op}{{{\rm op}}}

\newcommand{\argument}{{{\,\cdot\,}}}

\newcommand{\XX}{{\mathfrak X}}

\newcommand{\UU}{{\mathfrak U}}

\newcommand{\bB}{{\mathscr B}}
\newcommand{\C}{{\mathscr C}}

\newcommand{\F}{{\mathscr F}}

\renewcommand{\L}{{\mathscr L}}

\newcommand{\V}{{\mathscr V}}

\newcommand{\Uu}{{\mathfrak u}}

\newcommand{\zz}{{\mathbb Z}}
\newcommand{\hh}{{\mathbb H}}

\newcommand{\aaa}{{\mathbb A}}
\newcommand{\ttt}{{\mathbb T}}
\newcommand{\nn}{{\mathbb N}}

\newcommand{\qq}{{\mathbb Q}}

\newcommand{\cc}{{\mathbb C}}

\renewcommand{\O}{{\mathscr O}}

\newcommand{\aA}{{\mathscr A}}
\newcommand{\cC}{{\mathscr C}}

\newcommand{\fF}{{\mathscr F}}
\newcommand{\lL}{{\mathscr L}}

\newcommand{\hH}{{\mathscr H}}

\newcommand{\rR}{{\mathscr R}}

\newcommand{\st}{\mathrel{\mid}}

\newcommand{\tot}{\mathop{\rm tot}}

\newcommand{\Def}{\mathop{\rm Def}\nolimits}

\newcommand{\gr}{{\mathop{\rm gr}\nolimits}}

\newcommand{\Mod}{\mathop{{\rm Mod}}\nolimits}

\newcommand{\rank}{\mathop{\rm rank}\nolimits}

\newcommand{\red}{{\mathop{\rm red}\nolimits}}

\newcommand{\ob}{\mathop{\rm ob}}

\newcommand{\sheafhom}{\mathop{\rm {\mathscr{H}\!\mit{om}}}\nolimits}

\newcommand{\GL}{\mathop{\rm GL}\nolimits}

\newcommand{\ol}{\overline}

\newcommand{\doublearrowstack}[2]%
                      {{{{\scriptstyle#1}\atop{\textstyle\longrightarrow}}\atop{{\textstyle\longrightarrow}\atop{\scriptstyle#2}}}}
\newcommand{\rightleftarrowstack}[2]%
                      {{{{\scriptstyle#1}\atop{\textstyle\longrightarrow}}\atop{{\textstyle\longleftarrow}\atop{\scriptstyle#2}}}}
\newcommand{\leftrightarrowstack}[2]%
                      {{{{\scriptstyle#1}\atop{\textstyle\longleftarrow}}\atop{{\textstyle\longrightarrow}\atop{\scriptstyle#2}}}}

\newtheorem{thm}{Theorem}[section]
\newtheorem{cor}[thm]{Corollary}
\newtheorem{lem}[thm]{Lemma}
\newtheorem{prop}[thm]{Proposition}

\theorembodyfont{\upshape}
\newtheorem{defn}[thm]{Definition}

\newtheorem{rmk}[thm]{Remark}

\newenvironment{pf}{\begin{trivlist}\item[]{\sc Proof.}}%
            {\nolinebreak $\Box$ \end{trivlist}}

            {\nolinebreak $\Box$ \end{trivlist}}

\newcommand{\com}{\bullet}

\newcommand{\Proj}{\mathop{\mathbb{P}\rm{roj}}\nolimits}
\newcommand{\Gm}{{{\mathbb G}_m}}

\DeclareMathOperator\id{id}

\newcommand{\rk}{\mathop{\rm rk}\nolimits}

\newcommand{\spec}{\mathop{\rm Spec}\nolimits}

\newcommand{\Sym}{\mathop{\rm Sym}\nolimits}
\newcommand{\Ext}{\mathop{\rm Ext}\nolimits}
\newcommand{\Hom}{\mathop{\rm Hom}\nolimits}
\newcommand{\RHom}{\mathop{\rm RHom}\nolimits}

\newcommand{\hoch}{H\!H}

\newcommand{\noprint}[1]{}

\newcommand{\unsure}[1]{{\footnotesize #1}}

\def\unsure{\noprint}

\newcommand{\git}{\!\sslash\!}

\author{K. Behrend, B. Noohi}
\title{Moduli of non-commutative polarized schemes}

\begin{document}
\sloppy

\maketitle

\begin{abstract}
{We construct, using geometric invariant theory, 
a quasi-projective Deligne-Mumford stack of stable graded algebras. We
also construct a derived enhancement, which classifies twisted bundles
of stable graded $A_\infty$-algebras.  The tangent complex of the
derived scheme is given by graded Hochschild cohomology, which we
relate to ordinary Hochschild cohomology.  We obtain a
version of Hilbert stability for non-commutative projective schemes. 
}
\end{abstract}

\tableofcontents

\section*{Introduction}
\addcontentsline{toc}{section}{Introduction}

All our graded algebras will be unital and associative, with finite
dimensional graded pieces. 

We study derived moduli of graded  algebras. In the first part of this
paper, we construct a
differential graded stack $X$, classifying graded algebras of a fixed
dimension $\vec d=(d_1,d_2,\ldots)$. The construction is as a stack
quotient of a vector bundle of curved differential graded Lie algebras
over a linear space, divided by an algebraic `gauge' group.  The
construction is infinite-dimensional, but equal to the projective
limit of its finite dimensional truncations.

For a graded algebra $A$,
representing the point $P$ of $X$, the tangent complex of $X$ at $P$
has (shifted) Hochschild cohomology of $A$, computed with homogeneous
cochains of degree $0$, for cohomology groups: 
$H^i(\ttt_X|_P)=\hoch^{i+1}_\gr(A)$.  In other words, the derived
deformation theory of a graded algebra is given by its (shifted) graded
Hochschild cohomology. 

In the second section, we study graded Hochschild
cohomology in some detail, and relate it to more familiar
invariants. We do this for algebras $A$  `coming from geometry', by
which we mean that there exists a $\cc$-linear Grothendieck category
$\C$, an object $\O$ of $\C$, and an autoequivalence $s$ of $\C$,
satisfying suitable hypotheses, such that
$A=\bigoplus_{n\geq0}\Hom_\C(\O,s^n\O)$. In the commutative case,
this essentially means that
$A=\bigoplus_{n\geq0}\Gamma\big(X,\O(n)\big)$, for a
projective scheme $\big(X,\O(1)\big)$. 

Our results can be understood as supporting the idea that (under
certain hypotheses), the derived deformation theory of the graded
algebra $A$ coincides with the derived deformation theory of the triple
$(\C,\O,s)$.

In the last part of the paper we define a notion of {\em stability
}for graded algebras and construct a (derived) separated
Deligne-Mumford stack $\tilde X^s$ classifying stable graded algebras
of fixed dimension vector.  Our notion of stability comes from
geometric invariant theory for the finite-dimensional truncations
$\tilde X_{\leq q}$ of $\tilde X$. We expect that for many interesting
dimension vectors the stack $\tilde X^s$ (or at least an interesting
substack) will be of finite type.

In the commutative case, our notion of stability coincides with the
classical notion of Hilbert stability. Thus our stack $\tilde X^s$
extends the classical stack of Hilbert stable projective varieties
into the non-commutative world. 

We speculate that the stack $\tilde X^s$ (or a suitable open substack)
is a moduli stack of non-commutative projective schemes in the sense
of Artin-Zhang~\cite{AZ}.

As an example, 3-dimensional quadratic Artin-Schelter regular
algebras are semi-stable, and generically stable~\cite{BehHwa}.

\subsubsection{Notation and Conventions}

We work over a field of characteristic zero, which we shall denote by
$\cc$.  Unless specified otherwise, a graded vector space will refer
to a $\zz$-graded vector space. A graded vector space is {\em locally
  finite}, if each graded piece is finite dimensional over $\cc$. 

All our graded algebras will be unital and associative, locally finite
and graded in non-negative degrees. The component in degree zero will be assumed to
be equal to the ground field $\cc$.
  Often we will replace such an algebra $A$ by its
graded ideal $A_{>0}$ of elements of positive degree ($A$ can be
recovered from $A_{>0}$ in a canonical way). 

Our algebraic stacks will have affine diagonal, but we do not require
the diagonal to be of finite type, in general. 

We follow the Bourbaki convention that set inclusion (proper or not)
is denoted by `$\subset$'.

\subsubsection{Acknowledgements}

We started this project at the 10th Lisbon Summer Lectures in Geometry,
2009.  We thank Gustavo Granja for the hospitality during our
visit. 

This work was supported by a grant from the Royal Society under the International
Exchange Scheme.

\section{The derived stack of graded algebras}\label{secone}

\subsubsection{Gerstenhaber bracket}

Let $V=\bigoplus_{n>0} V_n$ be a locally finite positively graded
vector space. For $p\geq0$, let
$$L^p=\Hom_\gr(V^{\otimes(p+1)},V)$$
be the vector space of $(p+1)$-ary multilinear operations on $V$,
which  preserve degree. We have
$$L^p=\prod_{n}\Hom\big((V^{\otimes(p+1)})_n,V_n\big)\,,$$
which is a product of finite dimensional vector spaces. Therefore, it is an
affine $\cc$-scheme.  If $V$ is finite 
dimensional, it is an affine $\cc$-scheme of finite type. 

On $L=\bigoplus_{p\geq0} L^p$ we introduce a (non-associative) product
$\circ:L^p\otimes L^q\to 
L^{p+q}$ by the formula
$$(\mu\circ\nu)(a_0,\ldots,a_{p+q})=\sum_{i=0}^{p}(-1)^{iq}
\mu(a_0,\ldots,\nu(a_i,\ldots,a_{i+q}),\ldots,a_{p+q})\,.$$
We antisymmetrize and obtain the {\em Gerstenhaber bracket}
$$[\mu,\nu]=\mu\circ\nu-(-1)^{pq}\nu\circ\mu\,.$$
The pair $(L,[\,,])$ is a graded Lie algebra (see~\cite{gerstenhaber}).
It is finite dimensional, if $V$ is finite dimensional.

\paragraph{Augmentation.}
Sometimes it will be convenient to augment $L$ by putting a copy of
$\cc$ in degree $-1$, i.e., setting $L^{-1}=\cc$, and defining the
differential $L^{-1}\to L^{0}$ to be the map $\cc\to\Hom_\gr(V,V)$
given by the tautological graded endomorphism $\gamma$ of $V$, which is
multiplication by the degree, so $\gamma(\mu)=\deg(\mu)\mu$, for
homogeneous elements $\mu\in V$.  Define the bracket by
$[L^{-1},L]=0$.  The fact that the augmented object is a differential
graded Lie algebra follows from the fact that the tautological
endomorphism $\gamma$ is central in $L$. (This kind of
construction would not work with the identity in place of $\gamma$,
as the identity is not central.) We denote by $\tilde L$ the graded
Lie algebra obtained by dividing $L$ in degree 0 by the ideal
$\cc\gamma$.  Note that $\tilde L$ is quasi-isomorphic to the augmented
$L$.

\subsubsection{Maurer-Cartan equation}

The Maurer-Cartan equation for $L$ is 
$$\mu\circ\mu=0\,,\qquad \mu\in L^1\,.$$
Thus a Maurer-Cartan element  $\mu$ is a degree preserving binary
operation $\mu:V\otimes V\to V$, 
satisfying the equation
$$(\mu\circ\mu)(a,b,c)=0\,,\qquad\text{for all $a,b,c\in V$}\,.$$
Equivalently,
$$\mu(\mu(a,b),c)=\mu(a,\mu(b,c))\,,\qquad\text{for all $a,b,c\in
  V$}\,,$$
i.e., $\mu$ is associative.  Thus the Maurer-Cartan locus $MC(L)$ of
$L$ is the scheme of all associative graded products on $V$. It is a
closed subscheme of the affine scheme $L^1$. 

\subsubsection{Gauge group---Moduli stack}

The {\em gauge group }of $L$ is $G=\prod_nGL(V_n)$.  It is an affine
group scheme over $\cc$, and it is algebraic, if $V$ is finite
dimensional. Its Lie algebra is 
$L^0=\Hom_\gr(V,V)$. The gauge group acts on $L$ by conjugation, preserving the
Gerstenhaber bracket, and hence the Maurer-Cartan locus.  The {\em
  moduli stack }of $L$ is the stack quotient
$$X=[MC(L)/G]\,.$$ It classifies graded
associative products on $V$ up to change of basis in $V$.  In other
words, $X$ classifies graded associative algebras (without
unit), whose underlying graded vector space is isomorphic to $V$. The
stack $X$ is an algebraic stack with affine diagonal, although the
diagonal is not of finite type, unless $V$ is finite dimensional.

Let $d_i=\dim V_i$, and $\vec d=(d_1,d_2,\ldots)$.  The groupoid
$X(T)$, for a scheme $T$, is the category of bundles of graded
algebras of rank $\vec d$, parametrized by $T$. Such a bundle of
algebras is given by a graded vector bundle $\V=\bigoplus_{n>0}\V_n$ over
$T$, where $\rank \V_n=d_n$, endowed with $\O_T$-bilinear operations
$\V_i\otimes \V_j\to \V_{i+j}$, satisfying associativity.  (We can
always add to such a bundle of graded algebras a copy of $\O_T$ in
degree $0$, and make it into a bundle of unital algebras, in a
canonical way.)

\begin{defn}
For a graded algebra, we call the automorphisms $\phi_\lambda$, for $\lambda\in\Gm$,
given by $\phi_\lambda(a)=\lambda^{\deg a}a$ on homogeneous elements,
{\bf tautological}. The tautological automorphisms define the tautological
one-parameter group of automorphisms. 
\end{defn}

This leads to a modified moduli problem:
Denote by $\Gamma$ the central one-parameter subgroup of $G$ which
acts with weight $n$ on $V_n$, for all $n$. Let $\tilde G$ be
$G/\Gamma$.  The Lie algebra of $\tilde G$ is $\tilde L^0$. The group
$\tilde G$ acts by conjugation on $\tilde L$, and we call $\tilde G$
the gauge group of $\tilde L$. It is an affine group scheme over
$\cc$. Consider the quotient stack 
$$\tilde
X=[MC(L)/\tilde G]\,,$$ which is again an algebraic stack, the moduli
stack of $\tilde L$. 

We have a morphism of stacks $X\to \tilde X$, which is a
$\Gm$-gerbe. 

The moduli problem solved by $\tilde X$ is the following:
for a scheme $T$, the groupoid $\tilde X(T)$ is the groupoid of pairs
$(\XX,\V)$, where $\XX$ is a $\Gm$-gerbe over $T$, and
$\V=\bigoplus_{n>0}\V_n$ is an $\XX$-twisted vector bundle
on $T$, where $\V_n$ is $n$-twisted, and $\rank(\V_n)=d_n$, for all
$n>0$. Moreover, $\V$ is endowed with the structure of graded
algebra. We call such pairs $(\XX,\V)$ {\em twisted bundles of graded
  algebras}. 

For a review of twisted sheaves, see~\cite{twisted}.  Our terminology
is as follows: a quasi-coherent
$\XX$-twisted sheaf $\fF$ on $T$ is a quasi-coherent sheaf on $\XX$.
Such a sheaf decomposes naturally into a direct sum
$\fF=\bigoplus_{n\in\zz}\fF_n$, where on $\fF_n$ the natural inertia
action is equal to the $n$-th power of the linear action given by the
$\O_\XX$-module structure on $\fF$. If $\fF=\fF_n$, we refer to $\fF$
as $n$-twisted. If all $\fF_n$ are vector bundles, we call $\fF$ a
twisted vector bundle. 

If the components of the dimension vector $\vec d$ of $V$ are {\em strongly
  coprime}, by which we mean that there exists a $k$ such that 
$(d_1,2d_2,\ldots,kd_k)=1$, then the  gerbe $X\to\tilde X$ is
trivial. In this case, the universal twisted bundle of algebras can be
represented by a bundle of algebras. 
It
can be constructed by twisting the given action of $G$ on $V$ by the
character $\chi:G\to\Gm$, defined by
$\chi(g_1,g_2,\ldots)=\det(g_1)^{r_1}\det(g_2)^{r_2}\ldots$, where the
$r_i$ are such that $\sum i r_i d_i=1$.  The point is that this
twist does not affect the action on $L$, but it changes the action on
$V$ in such a way that it factors through $\tilde G$. 

\subsubsection{Derived moduli stack of algebras}

One of the simplest kinds of derived moduli stacks is given by a {\em
  bundle of curved differential graded Lie algebras } on a smooth
  algebraic stack (see~\cite{BCFHR}, for the definitions).  In the present case, the
  construction is as follows.

We start with the affine scheme $L^1$, and construct over it a bundle
of curved differential graded Lie algebras: the underlying graded vector
bundle $\lL=\bigoplus_{p\geq2}\lL^p$ is the trivial bundle over $L^1$, with fibre $L^p$ in degree
$p$, for $p\geq2$. The curvature map $L^1\to L^2$, given by $x\mapsto
x\circ x$, gives rise to a global section $f$ of $\lL^2$, the {\em
  curvature }of our bundle of curved differential graded Lie
algebras.  The twisted differential $d^\mu:\lL^i\to\lL^{i+1}$ is given
by $d^\mu=[\mu,\argument]$, in the fibre over $\mu\in L^1$.  The Lie
bracket on $\lL$ is constant over $L^1$, induced from the Gerstenhaber
bracket in each fibre. 

Then we notice that the gauge group action on $L^1$ lifts to an action
on all of $\lL$, preserving the structure of bundle of curved differential
graded Lie algebras.  Thus, this structure descends to the quotient
stack $M=[L^1/G]$, giving rise to a bundle of curved differential
graded Lie algebras over $M$.  From now on, let us reserve the notation
$(\lL,f,d^\mu,[\,,])$ for the descendant bundle on $M$. (If $V$ is
finite dimensional, each $\lL^p$  is a 
bundle of finite rank.) 

Our moduli stack $X$ is now realized as the closed substack $X\subset
M$, cut out by the vanishing of the curvature $f$ of
$\lL$.

In \cite{BCFHR}, it is explained how  a 
bundle of curved differential graded Lie algebras $(M,\lL)$ gives rise to a
differential graded stack, which we shall denote by $(M,\rR_M)$. In
fact, the curved differential graded Lie algebra structure on $\lL$
defines a differential graded co-algebra structure on $\Sym \lL[1]$,
which dualizes to a differential graded algebra structure on
$\rR_M=(\Sym \lL[1])^\vee$.   

We also get  a functor on differential
graded schemes: if $(T,\rR_T)$ is a differential graded scheme, 
we associate to it the set of  pairs
$(\V,\mu)$, where $\V$ is a graded vector bundle of dimension $\vec d$
over $T$, and $\mu$ is a global Maurer-Cartan element of the sheaf of
differential graded Lie algebras
$$\sheafhom_{\O_T}(\V^{\otimes\bullet+1},\V)\otimes_{\O_T}\rR_T\,.$$
This is the same thing as the
structure of a graded $A_\infty$-algebra 
on $\V\otimes_{\O_T}\rR_T$. 

We also have a bundle of curved differential graded Lie algebras over
$\tilde M=[L^1/\tilde G]$, giving rise to a differential graded stack
$(\tilde M,\rR_{\tilde M})$, whose underlying classical stack is
  $\tilde X$. This gives rise to the derived stack of twisted bundles
of graded $A_\infty$-algebras.

\subsubsection{Hochschild cohomology---Deformation theory}

Let us consider a point of $X$, represented by the Maurer-Cartan
element $\mu\in L^1$.  The derived stack $(M,\lL)$ gives rise to a complex of
vector bundles on $X$, the {\em tangent complex}, which governs
deformations and obstructions of morphisms from square zero extensions
of differential graded schemes. (For details, see~\cite{BCFHR}). At
the point $\mu$, this complex is our original graded Lie
algebra $L$, endowed with the twisted differential
$d^\mu=[\mu,\argument]$. This differential is the {\em Hochschild
  differential }of the associative algebra $(V,\mu)$. It makes
$(L,d^\mu,[\,,])$ into a differential graded 
Lie algebra.  

Explicitly, for $\alpha\in L^p$, $\alpha:V^{\otimes p+1}\to V$, we
have
\begin{multline*}
(d^\mu\alpha)(a_0,\ldots,a_{p+1})=\\
\alpha(a_0,\ldots,a_p)\,a_{p+1}+(-1)^pa_0\,\alpha(a_1,\ldots,a_{p+1})\\
-(-1)^p\sum_{i=0}^p(-1)^i\alpha(\ldots,a_ia_{i+1},\ldots)\,,
\end{multline*}
where we have written $\mu$ as concatenation.

The cohomology spaces
$$H^p(L,d^\mu)=\hoch^{p+1}_\gr(V,\mu)$$
are the {\em graded }Hochschild cohomology spaces of the graded
associative algebra 
$(V,\mu)$, computed with Hochschild cochains which are homogeneous of
degree zero.  Graded deformations/obstructions of the
graded algebra $(V,\mu)$  are given by $H^1(L,d^\mu)$ and
$H^2(L,d^\mu)$, respectively. 

Explicitly, if $\alpha:V^{\otimes 2}\to V$ is a 1-cocycle with
respect to $d^\mu$ (a Hochschild 2-cocycle), then 
\begin{equation}\label{coc}
\alpha(a,b)\,c-\alpha(a,bc)+\alpha(ab,c)-a\,\alpha(b,c)=0\,.
\end{equation}

The corresponding infinitesimal deformation of $(V,\mu)$ is given by
$V_{\epsilon}=V\oplus \epsilon V$ with multiplication $\ast$, which is
determined on $V\subset V_{\epsilon}$ by
$$a\ast b = ab + \epsilon\, \alpha(a,b)\,.$$
Associativity of $\ast$ follows from the cocycle
condition~(\ref{coc}). 

If $\beta:V\to V$ is a 0-cochain (a Hochschild 1-cochain), then
$\id+\beta:V_\epsilon\to V_{\epsilon}$ defines an isomorphism from
$\ast_{\alpha}$ to $\ast_{\alpha+d^\mu\beta}$. 

The infinitesimal deformation given by $\alpha$ extends to
$\cc[\epsilon]/\epsilon^3$, if and only if the {\em primary
  obstruction }$\alpha\circ\alpha$ vanishes in
$H^2(L,d^\mu)=HH^3_\gr(V,\mu)$. 

\subsubsection{Truncation}

Let $V\to V_{\leq q}$ be the truncation of $V$ into degrees less
than or equal to $q$, for a positive integer $q$. We will always
consider $V_{\leq q}$ as a quotient of $V$. Repeating the above
constructions with $V$ replaced by $V_{\leq q}$, we obtain a finite
dimensional graded Lie algebra $L_{\leq q}=\bigoplus L_{\leq q}^p$,
together with an epimorphism of graded Lie algebras $L\to L_{\leq
  q}$.
Let $X_{\leq q}$ and $\tilde X_{\leq q}$ denote the corresponding moduli
stacks, which are algebraic stacks of finite type, whose diagonal is
affine of finite type. We have
\begin{equation}\label{trunc}
X=\varprojlim_q X_{\leq q}\,,\qquad \tilde X=\varprojlim_q \tilde
X_{\leq q}\,.
\end{equation}

Let us write $M=[L^1/G]$ and $\tilde M=[L^1/\tilde G]$, etc.  Then we
have also $\tau_q:M\to M_{\leq q}$, and a morphism of bundles of
curved differential graded Lie algebras
$$\L\to \tau_q^\ast \L_{\leq q}\,,$$
for every $q$. 
Then
\begin{equation}\label{obvs}
\L=\varprojlim_q \tau_q^\ast\L_{\leq q}\,,
\end{equation}
(and a similar fact with tildes), as bundles of curved differential
graded Lie algebras. 

To state the compatibility with truncations on the level of
deformation theory, let $A=(V,\mu)$ be an algebra giving rise to a
point of $X$. Then we have
$$H^p(L,d^\mu)=\varprojlim_q H^p(L_{\leq q},d^\mu)\,,$$
as a direct consequence of (\ref{obvs}).  Thus, we also have
$$\hoch_\gr^p (A,A)=\varprojlim_q \hoch_\gr^p(A_{\leq q},A_{\leq q})\,.$$

\begin{rmk}
The projective system $L_q$ is a projective system of $\cc$-vector
spaces, and all transition maps are obviously surjective. The reason
to insist that we think of $V\to V_{\leq q}$ as a quotient map is
only to prove that the Hochschild boundary commutes with the maps of
the projective system.  A simple argument with $\varprojlim^1$ proves
that taking cohomology commutes with the projective limit.
\end{rmk}

\section{Graded Hochschild cohomology}

In this section we will relate graded Hochschild cohomology to more
familiar invariants.  We will do this for certain graded rings $S$ which `come
from geometry'.  By this we mean that $S$ is the homogeneous
coordinate ring of a `sufficiently amply polarized'  non-commutative
projective scheme $(\cC,A,s)$ in the sense of \cite{AZ}. Since our
hypotheses are going to diverge slightly from \cite{AZ}, we will call our
triples $(\cC,A,s)$  {\em
  polarized Grothendieck categories}, rather than non-commutative
  projective schemes.

We will define the concept of {\em reduced Hochschild cohomology }for
such a triple.  We apologize for this abuse of established
terminology.

\subsection{Hochschild cohomology of a polarized Grothendieck category}

\subsubsection{Preliminaries}

We summarize a result from \cite{LV}, which allows us to write down a
relatively small complex which computes the Hochschild cohomology of a
Grothendieck category.

Let $\cC$ be a $\cc$-linear Grothendieck category, and
$A:\Uu\to\cC$ a $\cc$-linear functor from a $\cc$-linear category
$\Uu$. This situation gives rise to the 
Yoneda functor 
$\cC\to\Mod(\Uu)$, where $\Mod(\Uu)$ is the category of right
$\Uu$-modules, i.e., the category of $\cc$-linear
functors
$\Uu^\op\to (\text{$\cc$-vector spaces})$.

We will need to assume that $\cC\to\Mod(\Uu)$ is fully faithful and
has an exact left adjoint. By the 
Gabriel-Popescu theorem, for this it suffices that 
$\{A(u)\}_{u\in  \ob\Uu}$ is a generating family for $\cC$, and
that $A:\Uu\to\cC$ is fully faithful.  By Theorem~1.2 of
\cite{LowenGP}, this latter 
condition can be weakened to 
\begin{items}
\item $A:\Uu\to\cC$ is faithful,
\item \label{local} for objects $u$, $v$ in $\Uu$, and a morphism  $f:A(u)\to A(v)$ in
  $\cC$, 
  there always exists a family of morphisms $u_i\to u$ in $\Uu$, such that
  $\coprod_i A(u_i)\to A(u)$ is an epimorphism in $\cC$, and $f|_{A(u_i)}\in
  \Uu$, for all $i$.
\end{items}

There exist (Section~1.10 in~\cite{Tohoku}) functorial injective
resolutions for the objects of 
$\Uu$. This 
means we have a 2-commutative diagram 
\begin{equation}\label{fir}
\vcenter{\xymatrix{
\Uu\rto^A \drto_E\drtwocell\omit{<-2>} & \cC\dto\\
& C^\com(\cC)}}\end{equation}
where  $C^\com(\cC)$ denotes the differential graded category of
complexes in $\cC$.
For every $u\in\Uu$, the homomorphism of
complexes $A(u)\to E(u)$ (given by the natural transformation
`$\Rightarrow$' in the diagram) is an injective resolution.

Denote by $\tilde E$ the $\Uu$-bimodule defined by the
functor
$E:\Uu\to C^\com(\cC)$.  We have
$$\tilde E(u,v)=\Hom^\com_\cC\big(E(u),E(v)\big)=
\RHom_\cC\big(A(u),A(v)\big)\,.$$  
We shall consider the Hochschild cochain complex
$C^\com(\Uu,\tilde E)$, see \cite[(2.4)]{LV}. It is the product total
complex of the double complex whose $p$-th column is given by 
\begin{multline*}
\prod_{u_0,\ldots, u_p}\Hom_\cc\bigg(
\Hom_\Uu(u_{p-1},u_{p})\otimes\ldots\otimes
\Hom_\Uu(u_0,u_1),\\
\Hom^\com_\cC\big(E(u_0),E(u_p)\big)\bigg)\,.\end{multline*}

\begin{prop}[\cite{LV}, Lemma 5.4.2]
The complex 
$C^\com(\Uu,\tilde E)$ computes the Hochschild cohomology of $\cC$ as
abelian category, and therefore governs the deformation theory of $\cC$
as abelian category.
\end{prop}

We will apply this result in the situation where $\{A(-n)\}_{n\in\nn}$
is a family of objects of $\cC$, such that for every $N$, the family
$\{A(-n)\}_{n<N}$ generates $\cC$.  We let $\Uu$ be the 
category whose objects are the negative
integers, and whose morphisms are given by 
$$\Uu(-m,-n)=\begin{cases}
\Hom_\cC\big(A(-m),A(-n)\big) & \text{if $-m\leq -n$}\\
0 & \text{if $-m>-n$}\end{cases}\,.$$
By construction, $\Uu$ comes with a faithful (although not necessarily
full) functor $A:\Uu\to\cC$, which satisfies Condition~\ref{local},
above.  The Hochschild complex $C^\com(\Uu,\tilde E)$ is given by 
\begin{multline}
\prod_{-n_0\leq\ldots\leq-n_p}\Hom_\cc\bigg(
\Hom\big(A(-n_{p-1}),A(-n_{p})\big)\otimes\ldots\otimes
\Hom\big(A(-n_0),A(-n_1)\big),\\ \label{column}
\Hom^\com_\cC\big(E(-n_0),E(-n_p)\big)\bigg)\,.\end{multline}

\subsubsection{Polarized Grothendieck categories}

Let $\cC$ be a $\cc$-linear Grothendieck category.  A {\bf polarization
}of $\cC$ is a pair $(s,A)$, where $s$ is an auto-equivalence of $\cC$,
and $A$ is an object of $\cC$, such that
\begin{items}
\item for every $N$, the family $\big(A(n)\big)_{n<N}$ generates $\cC$,
\item $\Ext^i_\cC\big(A,A(n)\big)=0$, if $n>0$, and $i>0$,
\end{items}
where we have written $s^nA=A(n)$.

In addition, we will make the assumption that $\Hom_\cC(A,A)=\cc$. 

For example, the Grothendieck category of quasi-coherent
$\O_X$-modules on a projective $\cc$-scheme $X$ is polarized by
$(\fF\mapsto\fF(1), \O_X)$, if $\O_X(1)$ is `sufficiently
ample'. It satisfies the additional assumption, if $X$ is connected. 

For another example, if $(\cC,A,s)$ is a finite-dimensional
non-commutative projective scheme in the
sense of \cite{AZ}, by which we mean that it satisfies the conditions
(H1), (H2), (H3), (H4), and (H5) of [ibid.], and has finite
cohomological dimension, then $(s,A)$ is a polarization of $\cC$, if
we replace $s$ by a sufficiently large power.  

As explained in \cite{AZ}, Proposition~4.2., we may, and shall, assume
that $s$ is an automorphism of $\cC$, rather than an autoequivalence. 

We choose functorial injective resolutions for the objects $A(-n)$,
$n\in\nn$, and use the complex $C^\com(\Uu,\tilde E)$, defined as
above~(\ref{column}), to compute the Hochschild cohomology of $\cC$.

\subsubsection{Reduced Hochschild cohomology}

Let $\tilde E^\ast$ be the same $\Uu$-bimodule as $\tilde E$,  except that we set 
$\tilde E(-1,-n)$ equal to zero:
$$\tilde E^\ast(-m,-n)=\begin{cases}
0 & \text{if $-m=-1$,}\\
\tilde E(-m,-n) & \text{otherwise}\,.
\end{cases}$$
By the definition of $\Uu$, we have that
$\tilde E^\ast$ is a bi-submodule 
of $\tilde E$. Let $\ol E$ be the quotient bimodule
$$\xymatrix{
0\rto & \tilde E^\ast\rto & \tilde E\rto & \ol E\rto & 0}\,.$$

Again, by the definition of $\Uu$, we have
for all $p$ that the $p$-th column of 
$C^\com(\Uu,\ol E)$ is a single copy of $\tilde
E(-1,-1)=\Hom^\com_\cC\big(E(-1),E(-1)\big)$. The
Hochschild differential is trivial, and therefore $C^\com(\Uu,\ol E)$
is quasi-isomorphic to $\tilde E(-1,-1)=\RHom_\cC(A,A)$. 

We call the cohomology of $C^\com(\Uu,\tilde E^\ast)$ the {\bf reduced
}Hochschild cohomology of $\cC$, with respect to the base object
$A$, notation $\ol\hoch^\ast(\cC,A)$. 

There is a short exact sequence of complexes
\begin{equation}\label{in1}
\xymatrix{
0\rto & C^\com(\Uu,\tilde E^\ast)\rto & C^\com(\Uu,\tilde E)\rto & 
C^\com(\Uu,\ol E)\rto & 0}\,,
\end{equation}
which gives rise to a long exact sequence in cohomology
\begin{equation}\label{redca}
\xymatrix{\rto &
\ol\hoch^\ast(\cC,A)\rto &\hoch^\ast(\cC)\rto &
\Ext_\cC^\ast(A,A)\rto^-{+1}&}\,.
\end{equation}

\begin{rmk}
As $\hoch^\ast(\cC)$ governs deformations of the abelian category  $\cC$, and
$\Ext_\cC^\ast(A,A)$ governs deformations of the object $A$ within
$\cC$ (see~\cite{TV}), the sequence (\ref{redca}) suggests that $\ol\hoch^\ast(\cC,A)$
governs the deformations of the pair $(\cC,A)$.  This motivates our
terminology.  We apologize for the somewhat ad-hoc definition, which
is motivated by its convenience for what follows.
\end{rmk}

\subsubsection{Graded Hochschild cohomology}

Define the unital graded $\cc$-algebra
$$S=\bigoplus_{n\geq0}\Hom_\cC\big(A,A(n)\big)\,,$$
and the {\em graded } differential graded $S$-bimodule 
$$M^\com=\bigoplus_{n\in\zz}\Hom^\com_\cC\big(E,E(n)\big)\,.$$
The grading coming from the autoequivalence $s$ will be called the
{\em projective }grading and will be denoted using lower indices, in
contrast to the cohomological grading, which is indicated with
superscripts.

We have the Hochschild complex of $S$ with values in  $M^\com$
$$C^\com(S,M^\com)$$
and the subcomplex $$C^\com_\gr(S,M^\com)$$ of projective degree $0$
cochains.  These are the cochains 
which preserve the projective degree. 

\begin{lem}\label{slem}
We have a short exact sequence of complexes of $\cc$-vector spaces
\begin{equation}\label{in2}
\xymatrix{
0\rto & C^\com_\gr(S,M^\com)\rto & C^\com(\Uu,\tilde E)\rto^-{1-s^{-1}} &
C^\com(\Uu,\tilde E)\rto & 0}\,.
\end{equation}
\end{lem}
\begin{pf}
During this proof we will disregard the vertical degree (the
coefficient degree) and consider only the horizontal degree (the
Hochschild degree).  Thus, $C^p(\Uu,\tilde E)$ will denote the $p$-th
column (\ref{column}) of $C^\com(\Uu,\tilde E)$. The same applies to
$C^\com_\gr(S,M^\com)$. 

A $p$-cochain $\chi\in C^p_\gr(S,M^\com)$, is a family
$(\chi_{\ell_1,\ldots,\ell_p})_{\ell_1,\ldots,\ell_p\geq0}$, where
$$\chi_{\ell_1,\ldots,\ell_p}:S_{\ell_p}\otimes\ldots \otimes
S_{\ell_1}\longrightarrow M^\com_{\ell_1+\ldots+\ell_p}$$
is a multilinear map.  We associate to  $\chi$  the $p$-cochain $\psi\in
C^p(\Uu,\tilde E)$ given by the family
$(\psi_{n_0,\ldots,n_p})_{n_0\geq\ldots\geq n_p\geq1}$, where
\begin{multline*}
\psi_{n_0,\ldots,n_p}:\Hom_\cC\big(A(-n_{p-1}),A(-n_p)\big)\otimes\ldots\otimes\Hom_\cC\big(A(-n_0),A(-n_1)\big)\\
\longrightarrow
\Hom^\com_\cC\big(E(-n_0),E(-n_p)\big)
\end{multline*}
is the multilinear operation given by 
$$
\psi_{n_0,\ldots,n_p}(\alpha_p,\ldots,\alpha_1)=s^{-n_0}\chi_{n_{0}-n_1,\ldots,n_{p-1}-n_p}(
s^{n_{p-1}}\alpha_p,\ldots,s^{n_0}\alpha_1)$$
Sending $\chi$ to $\psi$ defines the injective arrow in~(\ref{in2}). 

The functor $s^{-1}$ restricts to a fully faithful functor
$s^{-1}:\Uu\to\Uu$, and defines an endomorphism of the diagram
(\ref{fir}), and so induces an endomorphism $s^{-1}$ of $C^\com(\Uu,\tilde E)$. 
Given a $p$-cochain $\psi\in C^p(\Uu,\tilde E)$, the $p$-cochain
$s^{-1}\psi\in C^p(\Uu,\tilde E)$ is given by 
$$(s^{-1}\psi)_{n_0,\ldots,n_p}(\alpha_p,\ldots\alpha_1)=s\big(\psi_{n_0+1,\ldots,n_p+1}(s^{-1}\alpha_p,\ldots,s^{-1}\alpha_1)\big)\,.$$
So the condition $(1-s^{-1})\psi=0$ is equivalent to 
$$s\big(\psi_{n_0+1,\ldots,n_p+1}(s^{-1}\alpha_p,\ldots,s^{-1}\alpha_1)\big)=\psi_{n_0\ldots,n_p}(\alpha_p,\ldots,\alpha_1)\,.$$
Such a $\psi$ is the image of $\chi\in C^p_\gr(S,M^\com)$, with 
$$
\chi_{\ell_1,\ldots,\ell_p}(\alpha_p,\ldots,\alpha_1)=s^{n_0}\big(
\psi_{n_0,\ldots,n_p}(s^{-n_{p-1}}\alpha_p,\ldots,s^{-n_0}\alpha_1)\big)\,,
$$
where, for $i=0,\ldots,p$, we have used the abbreviation
$n_i=n+\sum_{j>i}\ell_j$, for an arbitrary $n\geq1$. This proves that
(\ref{in2}) is exact in the middle. 

To prove that $(1-s^{-1})$ is surjective, note that given $\phi$, the
equation $\phi=(1-s^{-1})\psi$ is equivalent to $s^{-1}\psi=\psi-\phi$,
which is a recursive equation for the components of $\psi$ in terms of
those $\psi_{n_0,\ldots,n_p}$ with $n_p=1$. 
\end{pf}

Sequences (\ref{in1}) and (\ref{in2}) exhibit two subcomplexes of
$C^\com(\Uu,\tilde E)$. 
In the intersection of $C^\com_\gr(S,M^\com)$ and
$C^\com(\Uu,\tilde E^\ast)$ inside $C^\com(\Uu,\tilde E)$, there is
$C^\com_\gr(S,S_{>0})$, giving rise to the commutative diagram of
complexes with exact rows and columns
\begin{equation}\label{ab}
\vcenter{\xymatrix{ 
C^\com_\gr(S,S_{>0})\rto\dto &  C^\com(\Uu,\tilde E^\ast)\rto\dto &
R\dto^\beta\\
C^\com_\gr(S, M^\com)\rto\dto & C^\com(\Uu,\tilde E)\rto^{1-s^{-1}}\dto &
C^\com(\Uu,\tilde E)\\
Q\rto^\alpha & C^\com(\Uu,\ol E)}}
\end{equation}

\begin{lem}
Both $\alpha$ and $\beta$ are quasi-isomorphisms. 
\end{lem}
\begin{pf}
In fact, the two claims are equivalent, so let us prove the one for
$\alpha$. Let us start by noting that in $C_\gr^\com(S,M^\com)$, we
can replace $M^\com$ by
$$M_{\geq0}^\com=\bigoplus_{n\geq0}\Hom_\cC^\com\big(E,E(n)\big)\,.$$
Consider the monomorphism of $S$-bimodules $S_{>0}\to
M^\com_{\geq0}$.  By the second condition that we require of
polarizations, the quotient of $M^\com_{\geq0}$ modulo $S_{>0}$ is
quasi-isomorphic to the bimodule $M_0=\Hom^\com_\cC(E,E)$, which exists
entirely in projective degree 0. 
It follows that $Q$ is quasi-isomorphic to $C^\com_\gr(S,M_0)$. But
for every $p$, we have $C^p_\gr(S,M_0)=M_0$. It follows that
$C^\com_\gr(S,M_0)$ is, in fact, quasi-isomorphic to
$\Hom^\com_\cC(E,E)=R\Hom_\cC(A,A)$. The same is true for
$C^\com(\Uu,\ol E)$. 

We have used the fact that graded Hochschild cohomology of $S$ is invariant
under quasi-isomorphisms of the coefficient bimodule.  This can be
reduced to the case of Hochschild cohomology of the category $\Uu$ via
Lemma~\ref{slem}, which applies to any $\Uu$-bimodule. 
\end{pf}

\begin{cor}
There is a distinguished triangle of complexes of $\cc$-vector
spaces
$$\xymatrix{
C^\com_\gr(S_{>0},S_{>0})\rto & C^\com(\Uu,\tilde E^\ast)\rto^-{1-s^{-1}}
& C^\com(\Uu,\tilde E) \rto^-{+1} &}\,,$$
and hence 
a long exact sequence in cohomology
\begin{equation}\label{ncseq}
\xymatrix{\rto &
\hoch_\gr^\ast(S_{>0},S_{>0})\rto & \ol\hoch^\ast(\cC,A)\rto^-{1-s^{-1}} &
\hoch^\ast(\cC)\rto^-{+1} &}\,.\end{equation}
\end{cor}
\begin{pf}
This is where we use the connectedness assumption that
$\Hom_\cC(A,A)=\cc$.  By this assumption, the normalized graded Hochschild
complex of $S$ with values in $S_{>0}$ is $C^\com(S_{>0},S_{>0})$. 
\end{pf}

Thus Diagram~(\ref{ab}) gives rise to a diagram of long exact
sequences in cohomology:
\begin{equation}\label{abprime}
\vcenter{
\xymatrix{
\hoch^\ast_\gr(S_{>0},S_{>0})\rto\dto & \ol\hoch^\ast(\cC,A)\rto^-{1-s^{-1}}\dto &
\hoch^\ast(\cC)\rto^-{+1}\ar@{=}[d] &\\
\hoch^\ast_\gr(S,M^\com)\rto\dto & \hoch^\ast(\cC)\rto^-{1-s^{-1}}\dto &
\hoch^\ast(\cC)\rto^-{+1} &\\
\Ext^\ast_\cC(A,A)\ar@{=}[r]\dto^{+1} & \Ext^\ast_\cC(A,A)\dto^{+1}\\
&
}}\end{equation}

\subsubsection{Heuristic Remarks}

Unfortunately, this result about 
$$(L,d^\mu)[-1]=C^\com_\gr(S_{>0},S_{>0})\,,$$
with notation $S_{>0}=(V,\mu)$, is only about the tangent complex of
our derived stack as a
complex, disregarding the structure of deformation functor, i.e. the
$L_\infty$-structure. We would like to make a few heuristic remarks, which
may explain the provenance of Diagram~(\ref{abprime}).

There is an octahedron of deformation functors
$$\xymatrix{
&\Def_\cC(A)\dto\ar@{=}[r] & \Def_\cC(A)\dto\\
\Def_\cC(s)\rto\ar@{=}[d] &\Def(\cC,s,A)\rto\dto & \Def(\cC,A)\dto\\
\Def_\cC(s)\rto & \Def(\cC,s)\rto & \Def(\cC)}$$
Then there is an isomorphism $\Def(\cC)=\Def_\cC(s)[1]$, so we can
rewrite this as
$$\xymatrix{
\Def_\cC(A)\dto\ar@{=}[r] & \Def_\cC(A)\dto\\
\Def(\cC,s,A)\rto\dto & \Def(\cC,A)\dto\rto& \Def(\cC)\ar@{=}[d]\\
\Def(\cC,s)\rto & \Def(\cC)\rto& \Def(\cC)}$$
and as
$$\xymatrix{ 
\Def(\cC,s,A)\rto\dto & \Def(\cC,A)\dto\rto& \Def(\cC)\ar@{=}[d]\\
\Def(\cC,s)\rto\dto & \Def(\cC)\rto\dto& \Def(\cC)\\
\Def_\cC(A)[1]\ar@{=}[r] & \Def_\cC(A)[1] & }$$

We believe that this latter diagram is, in fact, (\ref{abprime}), and
this justifies our suspicion that $C^\com_\gr(S_{>0},S_{>0})[+1]$ governs
the deformation theory of the triple $(\cC,A,s)$. From
Section~\ref{secone}, we know that $C^\com_\gr(S_{>0},S_{>0})[+1]$
governs the deformation theory of (non-unital) graded algebras.  This
is consistent with the Artin-Zhang philosophy that graded algebras are
just triples $(\cC,A,s)$ in disguise.

\subsection{Relative Hochschild cohomology (commutative case)}

In the commutative case, we can interpret graded Hochschild
cohomology of a graded ring as {\em reduced equivariant }Hochschild
cohomology.  We will introduce this concept, and prove results
analogous to the non-commutative case.

\subsubsection{Relative Hochschild cohomology for schemes}

Let $X$ be a separated scheme and $X\to Y$ a separated morphism of
algebraic stacks. Consider the diagonal morphism
$$\Delta:X\longrightarrow X\times_Y X\,,$$
which is a closed immersion of schemes. As for any closed immersion of
schemes, the derived category object $L\Delta^\ast\Delta_\ast\O_X$
splits off
$\hH^0(L\Delta^\ast\Delta_\ast\O_X)=\Delta^\ast\Delta_\ast\O_X=\O_X$,
and we write $(L\Delta^\ast\Delta_\ast\O_X)^\red$ for the complement
$\tau_{<0}(L\Delta^\ast\Delta_\ast\O_X)$. 

For a sheaf of $\O_X$-modules $\F$, we define the {\bf relative }Hochschild
cohomology of $X$ over $Y$ with values in $\F$ to be
$$\hoch^\ast_Y(X,\F)=R\Hom(L\Delta^\ast\Delta_\ast\O_X,\F)\,.$$
We also call 
$$\widebar\hoch^\ast_Y(X,\F)=R\Hom\big((L\Delta^\ast\Delta_\ast\O_X)^\red,\F\big)$$  
the {\bf   reduced }Hochschild cohomology of $X$ over $Y$ with values
in $\F$. For $\F=\O_X$, we use the usual abbreviations
$$\hoch_Y^\ast(X)=\hoch_Y^\ast(X,\O_X)\,,\qquad
\widebar\hoch_Y^\ast(X)=\widebar\hoch_Y^\ast(X,\O_X)\,.$$
We have
\begin{align*}
\hoch^\ast_Y(X,\F)&=\widebar\hoch^\ast_Y(X,\F)\oplus H^\ast(X,\F)\,,\\
\hoch^\ast_Y(X)&=\widebar\hoch^\ast_Y(X)\oplus H^\ast(X,\O_X)\,.
\end{align*}

\subsubsection{Equivariant Hochschild cohomology}

If $G$ is a reductive algebraic group, $\pi:P\to X$ is a principal
$G$-bundle, and $X\to BG$ the associated classifying morphism, then we
write $\hoch_G^\ast$ for $\hoch_{BG}^\ast$, and $\widebar\hoch^\ast_G$
for $\widebar\hoch^\ast_{BG}$, and speak of {\em
  equivariant }(reduced) Hochschild cohomology. 

\begin{prop}\label{grhoch}
For any quasi-coherent sheaf of $\O_X$-modules $\fF$, 
There is a natural $G$-action on $\hoch^\ast(P,\pi^\ast\fF)$, and we
have canonical isomorphisms
$$\hoch^\ast_G(X,\fF)=\hoch^\ast(P,\pi^\ast\fF)^G\,,\qquad
\widebar\hoch^\ast_G(X,\fF)=\widebar\hoch^\ast(P,\pi^\ast\fF)^G\,.$$
In particular,
$$
\hoch^\ast_G(X)=\hoch^\ast(P)^G\,,\qquad
\widebar\hoch^\ast_G(X)=\widebar\hoch^\ast(P)^G\,.$$
\end{prop}
\begin{pf}
Consider the cartesian diagram
$$\xymatrix{P\rrto^{\Delta'}\dto^\pi && P\times P\dto^{\tilde\pi}\\
X\rrto^\Delta && X\times_{BG} X}$$
and write $\aA=\pi_\ast\O_P$, so that $P$ is the relative spectrum of
the $\O_X$-algebra $\aA$ over $X$. 
By flat base change, we have 
$$\pi^\ast
L\Delta^\ast\Delta_\ast\O_X=L{\Delta'}^\ast{\Delta'}_\ast\O_P\,,$$
and therefore
$$
R\Hom(\pi^\ast L\Delta^\ast\Delta_\ast\O_X,\pi^\ast\fF)= 
R\Hom(L{\Delta'}^\ast{\Delta'}_\ast\O_P,\pi^\ast\fF)\,,$$
and by adjunction
$$\hoch^\ast_G(X,\aA\otimes_{\O_X}\fF)=\hoch^\ast(P,\pi^\ast\fF)\,.$$
We have a $G$-action on $\aA$, and the invariants are
$\aA^G=\O_X$. We get an induced action on $\aA\otimes_{\O_X}\fF$ with
invariants $\fF$, and an induced action
on 
$\hoch_G^\ast(X,\aA\otimes_{\O_X}\fF)$ with invariants
$\hoch_G^\ast(X,\fF)$. 
This proves the claim for usual Hochschild cohomology. The proof goes
through also in the reduced case.
\end{pf}

\subsubsection{Relation to ordinary Hochschild cohomology}

We specialize to the case $G=\Gm$. 

\begin{prop}
Let $X$ be a separated scheme and $X\to B\Gm$ a morphism. Denote the
diagonal $X\to X\times X$ by $\Delta$, and the diagonal $X\to
X\times_{B\Gm}X$ by $\tilde \Delta$. Then in
$D(\O_X)$ there are distinguished triangles
\begin{equation}\label{dt27}
\xymatrix{
L\Delta^\ast\Delta_\ast\O_X\rto^-{(t-1)} &
L\Delta^\ast\Delta_\ast\O_X\rto &
L\tilde\Delta^\ast\tilde\Delta_\ast\O_X\rto^-{+1}&}\,,
\end{equation}
and 
$$\xymatrix{
L\Delta^\ast\Delta_\ast\O_X\rto^-{(t-1)} &
(L\Delta^\ast\Delta_\ast\O_X)^\red\rto &
(L\tilde\Delta^\ast\tilde\Delta_\ast\O_X)^\red\rto^-{+1}&}\,.$$
\end{prop}
\begin{pf}
Let $X$ be a scheme, and 
$$\xymatrix{
\Gm\times X\rto^-\iota & \tilde R\rto^\pi & R}$$
a central extension of groupoids over $X$. 
The example which will concern us is given by $\tilde R=X\times_{B\Gm} X$, and
$R=X\times X$. 
Denote the identity
sections of $\tilde R$ and $R$ by $\tilde\Delta$ and $\Delta$,
respectively, and assume that $\Delta$ is a closed immersion. Then
$\iota$ is a closed immersion, as it is a pullback of $\Delta$. Let us
denote the 
identity of $\Gm\times X$ by $e$, and 
let $t$ be the standard coordinate on $\Gm$. 

We
have a short exact sequence of sheaves of $\O$-modules on $\Gm\times X$
$$\xymatrix{
0\rto & \O_{\Gm\times X}\rto^-{(t-1)} & \O_{\Gm\times X}\rto & e_\ast
\O_X\rto & 0}\,.$$
Applying $\iota_\ast$ we get the short exact
sequence 
$$\xymatrix{
0\rto & \iota_\ast\O_{\Gm\times X}\rto^-{(t-1)} &
\iota_\ast\O_{\Gm\times X}\rto & \iota_\ast e_\ast
\O_X\rto & 0}\,.$$
We have $\iota_\ast\O_{\Gm\times X}=\pi^\ast\Delta_\ast\O_X$, and
$\iota_\ast e_\ast\O_X=\tilde\Delta^\ast\O_X$, by the cartesian
diagram
$$\xymatrix{
\Gm\times X\rto^\iota\dto & \tilde R\dto^\pi\\
X\urto^{\tilde\Delta}\rto^\Delta & R}$$
So we can rewrite our exact sequence as
$$\xymatrix{
0\rto & \pi^\ast\Delta_\ast\O_{X}\rto^-{(t-1)} &
\pi^\ast\Delta_\ast\O_{X}\rto & \tilde\Delta_\ast 
\O_X\rto & 0}\,.$$
Now we apply $L\tilde\Delta^\ast$ to this exact sequence of
$\O$-modules on $\tilde R$, to obtain the distinguished
triangle~(\ref{dt27}).
\end{pf}

\begin{cor}
There are long exact sequences
$$\xymatrix{\rto&
\hoch^\ast_\Gm(X)\rto & \hoch^\ast(X)\rto & \hoch^\ast(X)\rto^-{+1}
&}\,,$$
and
$$\xymatrix{\rto&
\widebar\hoch^\ast_\Gm(X)\rto & \widebar\hoch^\ast(X)\rto & \hoch^\ast(X)\rto^-{+1}
&}\,.$$
\end{cor}

\subsubsection{A lemma on Hochschild cohomology of quasi-affine schemes}

If $X$ is quasi-affine, say $X\subset V=\spec A$, we can apply the usual
tilde construction to the Hochschild  complex $C_\com(A)$ of $A$.
We obtain a complex of quasi-coherent sheaves $C_\com(A)^\sim|_X$ on $X$, whose
component in degree $p$ is the free $\O_X$-module
$$C_p(A)^\sim|_X=\O_X\otimes_\cc A^{\otimes p}\,.$$
Removing the degree $0$ part from $C_\com(A)$ gives the reduced
Hochschild complex $\ol C_\com(A)$, and the associated complex of
quasi-coherent sheaves $\ol C_\com(A)^\sim|_X$, which is obtained
from $C_\com(A)^\sim|_X$ by removing the component in degree $0$.

\begin{lem}
In the derived category of $X$, the complex $C_\com(A)^\sim|_X$
represents the object $L\Delta^\ast\Delta_\ast\O_X$, where
$\Delta:X\to X\times X$ is the absolute diagonal. Moreover, the
complex $\ol C_\com(A)^\sim|_X$ represents
$(L\Delta^\ast\Delta_\ast\O_X)^\red$. 
\end{lem}
\begin{pf}
This follows from \cite{Swan}, where it is proved that on a
quasi-projective scheme the complex of non-quasi-coherent sheaves
$\C_\com^X$, which sheafifies the Hochschild complex, represents the
derived category object $L\Delta^\ast\Delta_\ast\O_X$. 

Then we have a canonical quasi-isomorphism
$$\xymatrix@1{C_\com(A)^\sim|_X\rto^-\sim &\C_\com^X}\,,$$
because Hochschild homology commutes with localization.
\end{pf}

Now suppose $\F=\tilde M|_X$, for an $A$-module $M$. 

\begin{lem}\label{spec}
We have spectral sequences
\begin{align*}
E_2^{pq}&=\hoch^p\big(A,H^q(X,\F)\big)\Longrightarrow
\hoch^{p+q}(X,\F)\,,\\
E_2^{pq}&=\widebar\hoch^p\big(A,H^q(X,\F)\big)\Longrightarrow
\widebar\hoch^{p+q}(X,\F)\,.
\end{align*}
\end{lem}
\begin{pf}
By the previous lemma, the derived 
category object $R\sheafhom(L\Delta^\ast\Delta_\ast \O_X,\F)$ can
be represented by the  complex $C^\com(A,\F)$ of sheaves on $X$, whose
degree $p$ component is 
$$\sheafhom_{\O_X}\big(\O_X\otimes_\cc A^{\otimes
  p},\F\big)=\Hom_\cc(A^{\otimes p},\F)\,,$$ 
i.e., an infinite product of copies of $\F$. 
It follows that Hochschild cohomology of $X$ with values in $\F$ is
equal to hypercohomology
$$\hoch^\ast(X,\F)=\hh^\ast\big(X,C^\com(A,\F)\big)\,.$$
This hypercohomology can be computed using a finite affine \v Cech
cover $\UU$  of $X$, because an infinite product of quasi-coherent sheaves is
acyclic over an affine scheme (even though not quasi-coherent in
itself). Thus
\begin{align*}
\hh^\ast\big(X,C^\com(A,\F)\big)& =
\tot\Big( \check C^\com\big(\UU,C^\com(A,\F)\big)\Big)\\
&=\tot\Big( C^\com\big(A,\check C^\com(\UU,\F)\big)\Big)\,.
\end{align*}
We now consider the double complex. Computing cohomology in the \v
Cech direction gives us $C^\com\big(A,H^q(X,\F)\big)$, because
infinite products are exact in the category of $A$-modules. Next,
computing cohomology in the Hochschild direction gives us
$\hoch^p\big(A,H^q(X,\F)\big)$, by definition.  Thus the desired spectral
sequence is the standard $E_2$ spectral sequence of our double
complex.

The proof is the same for the reduced case.
\end{pf}

\subsubsection{Graded Hochschild cohomology}

Let $A$ be a locally finite commutative graded $\cc$-algebra, such
that 
\begin{items}
\item $A$ is graded in non-negative degrees: $A=A_{\geq0}$, 
\item $A$ is connected: $A_0=\cc$,
\item $A$ is generated in degree $1$.
\end{items}

Let $V=\spec A$, and $Y=V\setminus 0$, where $0\in V$ is the vertex
defined by the homogeneous maximal ideal $A_{>0}$. Moreover, let
$X=Y/\Gm=\Proj A$, and denote the quotient map by 
$\pi:Y\to X$. Assume further that
\begin{items}\setcounter{enumi}{3}
\item for all $n>0$, the homomorphism  $A_n\to
  \Gamma\big(X,\O_X(n)\big)$ is bijective,
\item for all $q>0$ and $n>0$, we have $H^q\big(X,\O_X(n)\big)=0$.
\end{items}
Let us remark that 
$$H^q(Y,\O_Y)=H^q(X,\pi_\ast\O_Y)=\bigoplus_n
H^q\big(X,\O_X(n)\big)\,.$$
For example, if $X$ is a connected projective scheme, and $\O_X(1)$ is a sufficiently ample
line bundle, then $A=\bigoplus_i \Gamma\big(X,\O_X(i)\big)$ satisfies
our assumptions.

\begin{thm}\label{hochast}
We have
$$\ol\hoch^\ast_{\Gm}(X)=\hoch^\ast_\gr(A_{>0},A_{>0})\,.$$
\end{thm}
\begin{pf}
By Proposition~\ref{grhoch}, we have
$$\widebar\hoch^\ast_{\Gm}(X)=\widebar\hoch^\ast_\gr(Y)\,.$$
We can then use Lemma~\ref{spec} to determine
$\widebar\hoch^\ast_\gr(Y)$. In fact, $\Gm$ acts on the relevant spectral
sequence, and we get an induced spectral sequence of invariants
\begin{equation}\label{ssp}
E_2^{pq}=\widebar\hoch^p_\gr\big(A,H^q(Y,\O_Y)\big)\Longrightarrow
\widebar\hoch^{p+q}_\gr(Y)\,.
\end{equation}
To deal with the $E_2$-term, notice that, passing to the normalized
complex, we have
$$\hoch^p\big(A,H^q(Y,\O_Y)\big)=\hoch^p\big(A_{>0},H^q(Y,\O_Y)\big)\,,$$
This implies also
$$\widebar\hoch^p_\gr\big(A,H^q(Y,\O_Y)\big)
=\widebar\hoch^p_\gr\big(A_{>0},H^q(Y,\O_Y)\big)\,.$$  
For $q>0$ and $p>0$, we have
$$\hoch^p_\gr\big(A_{>0},H^q(Y,\O_Y)\big)=0\,,$$
because there are no graded cochains in the
relevant degrees (and taking invariants commutes with computing
Hochschild cohomology). So the $E_2$-term of the spectral
sequence~(\ref{ssp}) is entirely contained in the row $q=0$. We deduce
that 
$$\widebar\hoch^\ast_\gr(Y)=\widebar\hoch^\ast_\gr\big(A_{>0},H^0(Y,\O_Y)\big)^\gr\,.$$
We have
$$C^p_\gr\big(A_{>0},H^0(Y,\O_Y)\big)=C^p_\gr(A_{>0},A_{>0})\,,$$
and we conclude that
$\widebar\hoch^\ast_\gr(Y)=\widebar\hoch^\ast_\gr(A_{>0},A_{>0})$. 
\end{pf}

\begin{rmk}
This
argument would fail for non-reduced Hochschild cohomology, because the
corresponding spectral sequence  would also contain the
non-vanishing $n=0$ column.  This is the reason for working with
reduced Hochschild cohomology.
In fact, for Hochschild cohomology, we have
$$\hoch^\ast_{\Gm}(X)=\hoch^\ast_\gr(A_{>0},A_{>0})\oplus H^\ast(X,\O_X)\,.$$
\end{rmk}

\begin{cor}
There is a long exact cohomology sequence
\begin{equation}\label{cseq}
\xymatrix{\rto&
\hoch^\ast_\gr(A_{>0},A_{>0})\rto & \widebar\hoch^\ast(X)\rto & \hoch^\ast(X)\rto^-{+1}
&}\,.\end{equation}
This sequence is also the sequence (\ref{ncseq}), for $S=A$.
\end{cor}

\subsection{The smooth case}

\subsubsection{Hochschild-Kostant-Rosenberg}

We return to the case of a separated scheme $X$, with a separated
morphism $X\to Y$ to an algebraic stack $Y$, and assume that $X\to Y$
is smooth, of relative dimension $d$. The usual proof of the
Hochschild-Kostant-Rosenberg theorem 
goes through
and gives
\begin{align*}
L\Delta^\ast\Delta_\ast\O_X&=\bigoplus_{j=0}^d\Omega^j_{X/Y}[j]\,,\\
(L\Delta^\ast\Delta_\ast\O_X)^\red&=\bigoplus_{j=1}^d\Omega^j_{X/Y}[j]\,.
\end{align*}

\begin{cor}
For relative Hochschild cohomology, we have
\begin{align*}
\hoch^q_Y(X)=\bigoplus_{j=0}^dH^{q-j}(X,\Lambda^j T_{X/Y})\,\\
\widebar\hoch^q_Y(X)=\bigoplus_{j=1}^dH^{q-j}(X,\Lambda^j T_{X/Y})\,.
\end{align*}
\end{cor}

In particular, consider the case $Y= B\Gm$, and $X$ smooth. 
The bundles $\Lambda^j T_{X/B\Gm}$ can be related to the $\Lambda^j T_X$
by considering the short exact sequence of vector bundles on $X$
$$
\xymatrix@1{
0\rto & \O_X\rto & T_{X/B\Gm}\rto & 
T_X\rto & 0}\,$$
(the Euler sequence), which 
induces, for 
every $j>0$, another short exact sequence of vector bundles 
$$\xymatrix@1{
0\rto & \Lambda^{j-1} T_X\rto & \Lambda^j T_{X/B\Gm}\rto & \Lambda^{j}
T_X\rto & 0}\,.$$

If $A$ is a graded ring as in Theorem~\ref{hochast}, and $X=\Proj A$
is smooth of dimension $d$, then for $q>0$ we have
$$\hoch^q_\gr(A_{>0},A_{>0})=\bigoplus_{j=1}^{d+1} H^{q-j}(X,\Lambda^j
T_{X/B\Gm})\,.$$

\subsubsection{Further considerations}

We consider the case that $\big(X,\O_X(1)\big)$ is a smooth projective connected
scheme of dimension $d$. The polarization $\O_X(1)$ gives rise to the morphism $X\to
B\Gm$. Assume that $\O_X(1)$ is sufficiently ample, so that the
hypotheses of Theorem~\ref{hochast} are satisfied. Let $A$ be the
homogeneous coordinate ring of $X$. Then $A$ defines a point of the
derived moduli scheme of algebras constructed in
Section~\ref{secone}. The tangent complex at $X$ of the derived scheme
is 
$$\hoch^\ast_\gr(A_{>0},A_{>0})[1]=\bigoplus_{j=1}^{d+1} H^{\ast}(X,\Lambda^j
T_{X/B\Gm})[1-j]\,.$$
Therefore, the virtual dimension of the derived scheme at the point
$X$ is
$$1-\chi(X,\O_X)=(-1)^{1+\dim X}p_a(X)\,,$$
i.e., the arithmetic genus up to sign.

In this case, the beginning of the long exact sequence (\ref{cseq}), or
(\ref{ncseq}), is a direct sum of long exact sequences
as in Figure~\ref{fig1},
\begin{figure}
$$\xymatrix@R=.8pc{
0\rto &H^0(X,\O_X)\rto & H^0(X,\tilde T)\rto & 
H^0(X,T_X) \ar `[dr]`[dl]`[ddlll]`[ddll][ddll]& \\
&&&&\\
& H^1(X,\O_X)\rto & H^1(X,\tilde T)\rto & 
H^1(X,T_X)\ar `[dr]`[dddl]`[ddddll][ddddll]&\\
0\rto & H^0(X,T_X)\rto & H^0(X,\Lambda^2\tilde T)\rto &
H^0(X,\Lambda^2T_X) \ar `d[dl]`[ddddlll]`[ddddll][ddddll]&\\
&&&\\
&&&\\
& H^2(X,\O_X)\rto &H^2(X,\tilde T)\rto & H^2(X,T_X)\rto &\\
& H^1(X,T_X)\rto & H^1(X,\Lambda^2\tilde T)\rto &H^1(X,\Lambda^2
T_X)\rto &\\
0 \rto & H^0(X,\Lambda^2 T_X)\rto & H^0(X,\Lambda^3 \tilde T)\rto &
H^0(X,\Lambda^3 T_X)\rto &
}$$
\caption{}\label{fig1}
\end{figure}
where we have written $\tilde T$ for $T_{X/B\Gm}$. The left column
contains $\hoch^\ast(X)[-1]$, the middle column
$\ol\hoch_{\Gm}^\ast(X)=\hoch^\ast_\gr(A_{>0},A_{>0})$, and the right
column $\ol\hoch^\ast(X)$. 

Thus, the infinitesimal non-commutative polarized automorphisms of $X$ are
given by 
$$\widebar\hoch^1_\Gm(X)=H^0(X,\tilde T)\,.$$
This is equal to the
classical, commutative infinitesimal automorphisms of the pair
$\big(X,\O_X(1)\big)$. It is an extension of the kernel of $H^0(X,T_X)\to
H^1(X,\O_X)$ by $\cc=H^0(X,\O_X)$. 

The infinitesimal non-commutative polarized deformations of $X$ are
given by $\widebar\hoch^2_\Gm(X)$. This splits
up into two direct summands
$$\widebar\hoch^2_\Gm(X)=H^1(X,\tilde T)\oplus
H^0(X,\Lambda^2\tilde T)\,.$$

There is the classical, commutative part $H^1(X,\tilde T)$. This is an
extension of the kernel of $H^1(X,T_X)\to H^2(X,\O_X)$, i.e., the
infinitesimal deformations of $X$ lifting to the pair
$\big(X,\O_X(1)\big)$, by the cokernel 
of $H^0(X,T_X)\to H^1(X,\O_X)$, i.e., the infinitesimal
deformations of $\O_X(1)$, modulo those that come from infinitesimal
automorphisms of $X$. 

Then there is the non-commutative part $H^0(X,\Lambda^2\tilde
T)$. This is an extension of the kernel of $H^0(X,\Lambda^2 T_X)\to
H^1(X,T_X)$ by $H^0(X,T_X)$. 
The subspace $H^0(X,T_X)$ corresponds to non-commutative deformations
of the graded sheaf of algebras $\bigoplus_n\O(n)$ coming from
automorphisms of $X$, via the `twisted coordinate ring construction'.
The quotient space consists non-commutative deformations of the
structure sheaf (given by $H^0(X,\Lambda^2 T_X$), which map to zero in
$H^1(X,T_X)$.

The infinitesimal non-commutative polarized obstructions of $X$ are
given by $\widebar\hoch^3_\Gm(X)$, and split up into three parts.  The
classical, commutative part $H^2(X,\tilde T)$, and two non-classical
parts $H^1(X,\Lambda^2\tilde T)$ and $H^0(X,\Lambda^3\tilde T)$. In
particular, they contain $H^0(X,\Lambda^2 T_X)$ as a subspace. 

\begin{rmk}
For the obstruction theory to be perfect at $X$, i.e., for the higher
obstructions to vanish, we could require
$$H^i(X,\Lambda^j T_X)=0\,,\qquad\text{for all $i+j\geq 3$}.$$
For $X$ a curve this is always true. This leads to the speculation
that there may be interesting moduli
spaces of  non-commutative polarized curves, which admit virtual
fundamental classes.

For surfaces, this would give the three conditions
$$H^1(X,\Lambda^2 T_X)=0\,,\qquad H^2(X,T_X)=0\,,\qquad
H^2(X,\Lambda^2 T_X)=0\,.$$
\end{rmk}

\section{Stability for  graded algebras}

In this section we study the geometric invariant theory quotient
associated to the action of $G$ on $L^1$ (notation from Section~\ref{secone}). 
Because of (\ref{trunc}) we restrict to the case of truncated
algebras.  Then both $L^1$ and $G$ are of finite type, and we are in
a classical geometric invariant theory context.

\subsection{The GIT problem}\label{git}

Here we construct quasi-projective moduli schemes of finite graded
stable algebras. We start by 
setting up a Geometric Invariant Theory problem.

Let $q$ be a positive integer, $\vec d=(d_1,\ldots,d_q)$ a vector of
positive integers, and  $$V=\bigoplus_{i=1}^q V_i$$ a finite-dimensional
graded vector space of dimension $\vec d$. 
Let  $G=\prod_{i=1}^q\GL(V_i)$. We write elements
of $V$ as $x=(x_1,\ldots,x_q)$ and elements of $G$ as
$g=(g_1,\ldots,g_q)$.

Let $R=\Hom_\gr(V^{\otimes 2},V)$, with elements written as
$\mu=(\mu_{ij})_{ij}$, where $\mu_{ij}:V_i\otimes V_j\to
V_{i+j}$. Note that $\mu_{ij}\not=0$ only if $i,j\geq 1$ and
$i+j\leq q$. 
Consider  the left action of $G$ on
$R$ by conjugation.  More precisely, for $g\in G$ and $\mu\in R$, we
have
$$(g\ast\mu)_{ij}=g_{i+j}\circ \mu_{ij}\circ( g_i^{-1}\otimes
g_j^{-1})\,.$$

\begin{rmk}
This is not a space of quiver representations. So we cannot directly
quote results for moduli of quiver representations.  Although similar
techniques do apply.
\end{rmk}

There are two canonical one-parameter subgroups of $G$.  The anti-diagonal
$\Delta^{-1}:\Gm\to G$ acts by scalar multiplication (i.e., by weight
1) on $R$, and hence destabilizes every point of $R$. The other,
$\Gamma:\Gm\to G$ given by $\Gamma(t)=(t,t^2,\ldots,t^q)$ acts
trivially on $R$, prompting us to pass from $G$ to $\tilde
G=G/\Gamma$. 

\begin{defn}
We call a vector of integers
$\theta=(\theta_1,\ldots,\theta_q)$ a
{\bf stability parameter} if 
\begin{items}
\item $$\sum_{i=1}^q\theta_id_i<0\,,$$
\item $$\sum_{i=1}^q i\theta_i d_i=0\,.$$
\end{items}
\end{defn}

Any stability parameter defines a character $\chi_\theta$ of $\tilde G$ by 
$$\chi_\theta(g_1,\ldots,g_q)=\prod_{i=1}^q \det(g_i)^{\theta_i}\,.$$
The second condition on $\theta$ says that $\theta$ factors through
$G\to\tilde G$, and the first condition implies that
$\langle\chi,\Delta^{-1}\rangle>0$. 

We then linearize the action of $\tilde G$ on $R$ by taking the
trivial line bundle on $R$, and lifting the action to $R\times\cc$ by the
formula $g\ast(\mu,t)=(g\ast\mu,\chi(g)^{-1}t)$. Then the GIT quotient
of $R$ by $\tilde G$ is
$$R\git\tilde G=\Proj\bigoplus_{n=0}^\infty \Gamma(R)^{\tilde
  G}_{\chi^n}\,,$$
where 
$$\Gamma(R)^{\tilde G}_{\chi^n}=\{f\in\Gamma(R)\st
f(g\ast\mu)=\chi(g)^nf(\mu)\}$$ 
are the twisted invariants.  Note that the condition
$\langle\chi,\Delta^{-1}\rangle>0$ implies that
$\bigoplus_{n}\Gamma(R)_{\chi^n}^{\tilde G}$ is
non-negatively graded. 

The GIT quotient is a projective scheme,
because the affine quotient $\spec\Gamma(R)^{\tilde G}$ is reduced to
a point.

Let $R^s\subset R^{ss}\subset R$ be the open subsets of stable and
semi-stable points in $R$, respectively. Then $[R^s/\tilde G]$ is a
separated Deligne-Mumford stack with quasi-projective coarse moduli space
$R^s\git\tilde G$. Moreover, $R^{ss}\git\tilde G=R\git\tilde G$ is a
projective scheme, containing $R^s\git\tilde G$ as an open
subscheme. If $R^s=R^{ss}$, then $[R^s/\tilde G]$ is a proper
Deligne-Mumford stack with projective coarse moduli space $R\git\tilde
G$.

When we need to specify the stability parameter, we call points of
$R^s$ ($R^{ss}$) {\em $\theta$-(semi)-stable}.

\subsection{The Hilbert-Mumford criterion}

We recall the Hilbert-Mumford criterion (see Proposition~2.5 in \cite{QuiverKing}):

\begin{prop}[Hilbert-Mumford numerical criterion]
The point $\mu\in R$ is (semi)-stable (with respect to the
linearization given by $\chi$) if and only if for every
non-trivial one-parameter subgroup $\lambda$ of $\tilde G$, such that
$\lim_{t\to0}\lambda(t)\ast\mu$ exists in $R$, we have
$\langle\chi,\lambda\rangle >0$ $(\geq0)$. 
\end{prop}

\begin{prop}\label{hm}
The point $\mu\in R$ is $\theta$-(semi)-stable, if and only if for all
descending filtrations $V=V^{(0)}\supset V^{(1)}\supset\ldots$,
compatible with the lower grading, and satisfying the
conditions
\begin{items}
\item For $n$ sufficiently large, $V^{(n)}=0$, but $V^{(1)}\not=0$,
\item $(V^{(k)})$ does not dominate the tautological filtration,
  where to dominate the tautological filtration means that
$V^{(k)}\supset V_{\geq k}$, for all $k$,
\item $\mu(V^{(i)},V^{(j)})\subset V^{(i+j)}$, for all $i$, $j$,
\end{items}
we have
$$\sum_{i=1}^q \theta_i w_i>0\quad \text{$(\geq0)$}\,.$$
Here $w_i=\sum_{m\geq1} \dim
  V_i^{(m)}$ is the {\bf weight function }of the filtration $V^{(k)}$.
\end{prop}
\begin{pf}
A one-parameter subgroup of $\tilde G$ is the same thing as a
one-parameter subgroup of $G$, up to translation by
$\Gamma$. One-parameter subgroups of $G$ are the same thing as
gradings on each of the $V_i$, which we denote by upper indices
$V_i=\bigoplus_m V_i^m$. The  upper grading $V=\bigoplus_{i,m}V^m_i$
gives rise to the same one-parameter subgroup of
$\tilde G$ as the upper grading $V=\bigoplus_{i,m}V^{m+i}_i$. Thus we
call the upper gradings $\bigoplus V_i^m$ and $\bigoplus V_i^{m+i}$
{\em equivalent}. The upper grading defined by $V=V^0$, as well as all
equivalent upper gradings are called {\em trivial}, as they correspond
to the trivial cocharacter of $\tilde G$. In each equivalence
class there is a unique upper grading such that no weights are
negative, but there exists a non-zero space $V^m_i$ with $m<i$. Let us
call such an upper grading {\em standard}. 

Now let $\mu\in R$ be given. Let $\lambda$ be a one-parameter subgroup
of $G$, corresponding to the double grading $V=\oplus V_i^m$ on
$V$. Then $\lim_{t\to 0}\lambda(t)\ast\mu$ exists in $R$, if and only
if none of the  $\lambda$-weights of $\mu$ are negative. This is
equivalent to $\mu$ preserving the descending filtration given by
$V^{\geq n}=\bigoplus_{m\geq n} V^m$, by which we mean that
$\mu(V^{\geq m}, V^{\geq n})\subset V^{\geq m+n}$. Note that this
condition is preserved under equivalence of upper gradings, even
though the upper filtration itself changes in the equivalence class. 

Now suppose that $\mu\in R$ preserves the filtration $V^{\geq n}$,
given by $\lambda$. Then 
$$\langle\chi,\lambda\rangle=\sum_{i=1}^q\theta_i\sum_{m}m\dim
V^m_i\,.$$
Note that for standard upper gradings, we have $V\subset
V^{\geq 0}$, and hence
$$\sum_m m\dim V^m_i=\sum_{m\geq1} \dim V^{\geq m}_i\,,$$
so that we have
$$\langle\chi,\lambda\rangle=\sum_{i=1}^q\theta_i\sum_{m\geq1}\dim
V_i^{\geq m}\,.$$
Thus we conclude that $\mu\in R$ is stable if and only if for every
descending filtration $V=V^{(0)}\supset V^{(1)}\supset\ldots$
(compatible with the lower grading), satisfying
\begin{items}
\item (non-trivial) $V^{(1)}\not=0$, but $V^{(n)}=0$, for $n\gg0$,
\item (standard) there exists a $k$, such that $V^{(k)}\not\supset
  V_{\geq k}$,
\item $\mu(V^{(m)},V^{(n)})\subset V^{(m+n)}$,
\end{items}
we have $\sum_{i=1}^q\theta_i \sum_{m\geq1} \dim
  V_i^{(m)}>0$ ($\geq0$).
\end{pf}

\subsection{Reformulation using test configurations}

Suppose now that $A$ is an associative and unital graded algebra, which
is a locally finite and connected, with $A_0=\cc$. 

\subsubsection{Test configurations for $A$}

\begin{defn}
A {\bf test configuration }for $A$ is 
a bundle of graded unital algebras $\bB$ (as defined in
  Section~\ref{secone})
over the affine line $\aaa^1$, together
with a 
$\Gm$-action on the bundle  $\bB$, lifting the natural action on
$\aaa^1$, such that the restriction of $\bB$ to 
$\Gm\subset\aaa^1$ is $\Gm$-equivariantly isomorphic to the constant
bundle with fibre $A$.

Two test configurations for $A$ are {\bf equivalent}, if one can be obtained
from the other by multiplying the $\Gm$-action by a suitable  power of
the tautological action. 
A test configuration for $A$ is {\bf trivial}, if it is equivalent to 
a $\Gm$-equivariantly constant test configuration.
\end{defn}

The special fibre $\bB|_0$ of a test configuration is a graded algebra
with the same Hilbert function as $A$, endowed with a
$\Gm$-action. The weight of the $\Gm$-action on the graded piece of
degree $i$ of $\bB|_0$ is denoted $w_i$, and the function 
$$F(i)=\frac{w_i}{i\dim A_i}\,,$$
defined for $i>0$, 
is the {\bf Futaki  function }of the test configuration $\bB$. (It
takes values in $\qq\cup\{\infty\}$.) The
Futaki functions of two equivalent test configurations differ by a
constant integer. The Futaki function of a trivial test configuration
is a constant integer.

A test configuration for $A$, together with a $\Gm$-equivariant
trivialization of its restriction to $\Gm\subset \aaa^1$,  is the same
thing as a doubly graded $\cc[t]$-subalgebra 
$$B=\bigoplus_{k\in\zz} B^{(k)}t^{-k}\subset
A[t,t^{-1}]\,,$$
such that every $B_i\subset A_i[t,t^{-1}]$ is a finitely generated
$\cc[t]$-submodule of rank $\dim A_i$.
The test configuration is trivial, if and only if there exists an
$\ell\in\zz$, such that 
$$B^{( k)}_i=\begin{cases}
A_i &\text{if $k\leq i\ell$}\\
0 & \text{if $k>i\ell$}\end{cases}$$
For our purposes it will not be important to distinguish between a
test configuration and one with $\Gm$-equivariant trivialization over
$\Gm\subset\aaa^1$, and so we will identify test configurations with
doubly graded $\cc[t]$-algebras $B\subset A[t,t^{-1}]$ such that $\rk
B_i=\dim A_i$, for all $i$. 

By definition, {\bf generators }of a test configuration $\bB$ for $A$ are
generators for the algebra of global sections
$B=\Gamma(\aaa^1,\bB)$ as $\cc[t]$-algebra. 

\begin{rmk}
If $A$ admits
a finitely generated test configuration, then $A$ is finitely generated,
itself.
\end{rmk}

\subsubsection{Admissible test configurations}

\begin{defn}
A test configuration is called {\bf admissible}, if it is equivalent
to a test configuration which can be written as
$$A[t]\subset B\subset A[t,t^{-1}]\,.$$
\end{defn}

Let us suppose $B$ is an admissible test configuration written in this
way. We have
$$w_i=\sum_{k>0}\dim B^{(k)}_i\,.$$
Moreover, 
\begin{items}
\item every $B^{(k)}$ for $k>0$ is a two-sided ideal in $A$, 
\item $A\supset B^{(1)}\supset B^{(2)}\supset\ldots$,
\item $B^{(k)} B^{(\ell)}\subset B^{(k+\ell)}$, for all $k,\ell>0$,
\item for every $i>0$, there exists an $\ell>0$, such that
  $B_i^{(k)}=0$, for all $k\geq\ell$.
\end{items}

\begin{defn}
We call a sequence of two-sided ideals $\big(I^{(k)}\big)_{k>0}$ in
$A$ satisfying these  conditions an {\bf admissible }family of
ideals in $A$.
\end{defn}

An admissible family of ideals $\big(I^{(k)}\big)_{k>0}$ defines a
test configuration by
$$B=\bigoplus_{k\in\zz} I^{(k)}\,t^{-k}\,,$$
where we set $I^{(k)}=A$, for all $k\leq0$. 
The special fibre of this test configuration is
$$B/tB=\bigoplus_{k\geq0}I^{(k)}/I^{(k+1)}\,.$$

\begin{rmk}
If $A$ is finitely generated, then every test configuration for $A$ is
admissible.
\end{rmk}

\subsubsection{Standard admissible test configurations}

If a test configuration is admissible, there is a unique equivalent
test configuration with the properties
\begin{items}
\item $A[t]\subset B$,
\item $\bigoplus_{k\in\zz} A_{\geq k}\,t^{-k} \subsetneqq B$.
\end{items}
Such a test configuration is called {\bf standard admissible}. 

A test configuration is standard admissible if and only if the
corresponding admissible family of ideals does not contain the
{\bf tautological }admissible family given by $I^{(k)}=A_{\geq k}$.  Such an admissible
family of ideals is called {\bf standard admissible}. 

\subsubsection{Stability}

Now let us return to the setup of \ref{git}. Suppose that $\mu\in R=L^1$
is a Maurer-Cartan element, so that $A=(V,\mu)$ is a graded algebra
with $A_i=0$, for $i>q$. 

\begin{prop}
The Maurer-Cartan element $\mu$ is $\theta$-(semi)-stable if and
only if, for every non-trivial test
configuration for $A$, the weights 
$w_i$ satisfy $$\sum_i\theta_i w_i >0\quad (\geq0)\,.$$
\end{prop}
\begin{pf}
As $A_i=0$ for $i\gg 0$, all test configurations for $A$ are
admissible. Because of $\sum_ii\theta_id_i=0$, the stability condition
$\sum_i\theta_i w_i>0$ $(\geq0)$ is independent of the choice of a
test configuration within its equivalence class.  So to test the condition of
this proposition it is sufficient to use standard admissible test
configurations. To conclude, we remark that non-trivial standard
admissible test configurations correspond exactly to the filtrations
of $V$, which are tested in Proposition~\ref{hm}.
\end{pf}

This proposition motivates the following definition.

\begin{defn}
A finite graded algebra $A$ is {\bf $\theta$-(semi)-stable},
if 
\begin{items}
\item $\sum_i \theta_i\dim A_i<0$,
\item $\sum_i i\theta_i\dim A_i=0$,
\item for every 
non-trivial test configuration for $A$, the weights
satisfy $$\sum_i\theta_i w_i >0\quad (\geq0)\,.$$ 
\end{items}
\end{defn}

\begin{prop}
To test (semi)-stability of $A$, it suffices to check standard
admissible families of ideals in $A$.
\end{prop}

\subsection{Standard stability parameters}

We fix a dimension vector $(d_1,\ldots,d_q)$ and a stability parameter
$\theta$, as above, and study $\theta$-stability of  finite graded
algebras $A$ with $\dim A_i=d_i$. 

Let us remark that there is no a priori reason to expect
complete moduli of stable algebras:

\begin{rmk}
We can eliminate $\theta_1$ from the stability condition. The
stability parameter condition becomes
$$\sum_{i=2}^q(i-1)d_i\theta_i>0\,,$$
and as stability condition we obtain
$$\sum_{i=2}^q(d_1w_i-id_iw_1)\theta_i>0\quad\text{($\geq0$)}\,,$$
or
$$\sum_{i=2}^q\big(F(i)-F(1)\big)id_i\theta_i>0\quad\text{($\geq0$)}\,.$$
We see that no matter the choice of stability parameter $\theta$, an
admissible sequence of ideals with constant Futaki function will
always violate stability. The Futaki function being constant means
that 
$$w_k=k d_k \frac{w_1}{d_1}\,,\qquad\text{for all $k\geq1$}\,.$$
There is no a priori reason why $d_1$ should not divide $w_1$, and so
there is no divisibility condition on the dimension vector
$(d_1,d_2,\ldots)$ which would exclude the possibility of strictly
semi-stable objects.  Therefore, there is no such condition that would
assure a projective coarse moduli space of stable algebras. 
\end{rmk}

For certain stability parameters, stability implies generated in
degree 1:

\begin{prop}\label{gen1}
Suppose that $\theta_1<0$ and $\theta_i\geq0$, for all $i>1$.  Then
$\theta$-stable algebras are generated in degree 1.  If, in addition,
$\theta_i>0$, for all $i>1$, then $\theta$-semi-stable algebras are
generated in degree 1.
\end{prop}
\begin{pf}
Write $I=A_{\geq1}$, and consider the admissible sequence of
ideals of powers of $I$, given by $I^{(k)}=I^k$, for all
$k\geq1$. Assume that {\em not }$I^k=A_{\geq k}$, for all $k$. Then
$(I^k)$ is properly contained in the tautological filtration, and
hence does not dominate it. Thus $(I^k)$ is standard admissible. 

If $A$ is $\theta$-stable, then $(-\theta_1)
w_1<\sum_{i>1}\theta_iw_i$. This implies $(-\theta_1)
d_1<\sum_{i>1}\theta_i w_i$, and hence
$\sum_{i>1}i\theta_id_i<\sum_{i>1}\theta_i w_i$. This is a
contradiction, because $w_i\leq id_i$, for all $i$. Thus $I^k=A_{\geq
  k}$, for all $k$, which implies that $A$ is generated in degree
1. 

To prove the additional claim, assume that $V$ is
$\theta$-semi-stable. Then we can still conclude that
$\sum_{i>1}i\theta_id_i\leq\sum_{i>1}\theta_iw_i$.  Thus, from
$w_i\leq id_i$, and the fact that none of the $\theta_i$ vanish, we
conclude that $w_i=i d_i$, for all $i>1$. Again, we reach a
contradiction, proving that $A$ is generated in degree 1.
\end{pf}

\begin{rmk}
We have, in both cases, proved that any admissible sequence of ideals which
is contained in the tautological one, and satisfies $I^{(1)}_1=A_1$,
is necessarily the tautological sequence. 
\end{rmk}

\begin{prop}\label{coq}
Suppose that  we have 
$\theta_i\leq0$, for all $i<q$. Then every 
$\theta$-stable algebra has no non-zero ideal $I$, which vanishes in
degree $q$. If $\theta_i<0$, for all $i<q$, we can reach the same
conclusion for $\theta$-semi-stable algebras. 
\end{prop}
\begin{pf}
In fact, if we assume that $I^{(k)}$ is an admissible sequence of
ideals which vanishes in degree $q$,  we can conclude that
$I^{(k)}=0$, for all $k\geq1$, under either of the two assumptions. 
\end{pf}

\begin{rmk}
If $\theta_q=0$, there are no stable algebras.
\end{rmk}

\begin{rmk}
If we want the assumptions of both Propositions~\ref{gen1}
and~\ref{coq} to hold, we need to have $\theta_1<0$, and
$\theta_q>0$, as well as
$\theta_i=0$, for all $1<i<q$.  For the conclusions to hold, we need
to assume stability, not just  semi-stability.
\end{rmk}

\begin{defn}
The stability parameter $\theta$ is {\bf standard}, if  $\theta_1$ and
$\theta_q$ are the only non-zero $\theta_i$.
\end{defn}

For a standard stability parameter $\theta$, the stability condition
reads
$$\theta_q w_q>(-\theta_1)w_1\quad (\geq)\,.$$
This is equivalent to
$$F(q)>F(1)\quad(\geq)\,,$$
which is  independent of the sizes of $\theta_1$ and $\theta_q$.

When not specified otherwise, we always work with a standard
stability condition, and make the following definition.

\begin{defn}
Let $A$ be a finite graded algebra, graded in the interval
$[0,q]$. Then $A$ is
called {\bf  (semi)-stable}, if for every non-trivial  
test configuration for $A$, the Futaki function satisfies $F(q)>F(1)$
$(\geq)$.  It suffices to check admissible families of ideals, or
standard admissible sequences of ideals. 
\end{defn}

\begin{cor}\label{coro}
Suppose $A$ is stable. Then $A$ is generated in degree $1$, and
has no non-trivial two-sided ideals which vanish in degree $q$.
\end{cor}

\subsubsection{Moduli}

Consider the dimension vector $\vec d=(d_1,\ldots,d_q)$, and the
associated stack of twisted bundles of graded algebras of
dimension $\vec d$, which we called $\tilde X_{\leq q}$ in
Section~\ref{secone}. Let $\tilde X^s_{\leq q}$ be the open substack of stable
algebras.  It is a closed substack of the quotient stack $[R^s/\tilde G]$, and
it is a separated Deligne-Mumford stack with quasi-projective coarse
moduli space, which is a closed subscheme of  $R^s\git \tilde G$. The
$\cc$-points of this coarse 
moduli space correspond in a one-to-one fashion to isomorphism classes
of stable algebras of dimension $\vec d$.

\subsection{Unbounded algebras}

For simplicity, we will only consider stability, not
semi-stability. In view of Corollary~\ref{coro}, we will  only
consider algebras generated in degree $1$.

\begin{prop}\label{eqcon}
Fix an integer $q>1$, and let $A$ be a  graded algebra,
finitely generated in degree $1$. The following 
are equivalent 
\begin{items}
\item 
For every test configuration $\bB$ for $A$, whose truncation  $\bB_{\leq
  q}$ is a non-trivial test configuration for $A_{\leq q}$,  the
Futaki function satisfies $F(q)>F(1)$.
\item 
For every  non-trivial test configuration for $A$ generated in degrees
$\leq q$, the Futaki
function satisfies $F(q)>F(1)$.
\item For every non-trivial test configuration for $A$ generated in
  degree $1$, the Futaki function satisfies $F(q)>F(1)$.
\item \label{filter}For every filtration  $A_1\supsetneq
V^{(1)}\supset\ldots\supset V^{(r)}\supsetneq 0 $ of  $A_1$
by vector subspaces, the admissible
sequence of ideals generated by 
$\{V^{(k)}\}$ has a Futaki function which satisfies $F(q)>F(1)$.
\item The truncation $A_{\leq q}$ is stable. 
\end{items}
\end{prop}
\begin{pf}
The fact that (i) implies (ii), follows because if a test
configuration $\bB$ for $A$ is generated in degrees $\leq q$, and is non-trivial,
then also its truncation $\bB_{\leq q}$ is non-trivial. 

Obviously, (ii) implies (iii). 

Next we claim that  (iii) implies (iv). Here we will use that $A$ is
generated in degree $1$. The  admissible sequence of
ideals generated by the  
filtration $\{V^{(k)}\}$ is the smallest admissible sequence of ideals
$\{I^{(k)}\}$, with $I^{(k)}_1=V_1^{(k)}$, for all $k>0$.
The corresponding test configuration $B$ is generated as
$\cc[t]$-algebra by $\bigoplus_{k\geq0} V^{(k)}t^{-k}\subset A_1[t,t^{-1}]$ inside
$A[t,t^{-1}]$, if we set $V^{(0)}=A_1$.  It is generated in degree
$1$. Thus, (iii) implies (iv). 

Now let us assume that  (iv) is satisfied. To prove (v), i.e., that
$A_{\leq q}$ is stable, it suffices to check all standard admissible
test configurations for $A_{\leq q}$. Among these, it suffices to
check those that are generated in degree $1$, because adding
generators in higher degree can only increase $F(q)$, without
affecting $F(1)$. But non-trivial standard admissible test
configurations generated in degree $1$ are all generated by a
filtration $\{V^{(k)}\}$ as in (iv).  

Finally, the fact that (v) implies (i) is, again, trivial.
\end{pf}

\begin{rmk}\label{further}
For a given dimension $d_1$ of $A_1$, in Condition~(iv), we can further reduce to
considering only flags whose dimensions $(\dim V^{(1)}, \dim
V^{(2)},\ldots)$ come from a finite list of integer sequences, but as
we currently have no use for this fact, we will not prove it here.
\end{rmk}

\begin{defn}\label{defstab}
A connected graded algebra, finitely generated in degree $1$,  is called {\bf
  $q$-stable}, if any of the equivalent conditions in
Proposition~\ref{eqcon} is satisfied. It is called {\bf stable}, if
there exists and $N>0$, such that it is $q$-stable for all $q\geq
N$. 
\end{defn}

\subsubsection{Commutative case}

Suppose that $(Y,\O_Y(1))$ is a connected projective $\cc$-scheme, such that
$H^i(Y,\O(1))=0$, for all $i>0$. Let $A$ be the homogeneous
coordinate ring of $(Y,\O_Y(1))$.  This is the image of $\Sym
\Gamma(Y,\O(1))\to \bigoplus_{i\geq0} \Gamma(Y,\O(n))$, and is a
connected graded algebra, generated in degree $1$.

\begin{prop}
The polarized scheme $(Y,\O_Y(1))$ is Hilbert stable if and
only if $A$  is stable in the sense of Definition~\ref{defstab}.
\end{prop}
\begin{pf}
For the definition of Hilbert stability (more precisely, Hilbert
stability with respect to $r=1$), see \cite{RT}.  By definition, the
Hilbert stability of $(Y,\O_Y(1))$ is tested against all 
filtrations of $A_1=\Gamma(Y,\O_Y(1))$, exactly as in
Proposition~\ref{eqcon}~\ref{filter}. This immediately implies the
result.
\end{pf}

\unsure{
\begin{lem}
If $B$ is admissible, and $A[t]\subset B\subset A[t,t^{-1}]$, then $B$ is generated in
degree $1$ as $\cc[t]$-algebra, if and only if the 
$\cc$-subalgebra $\bigoplus_{k>0}B^{(k)}t^{-k}\subset A[t,t^{-1}]$ is  generated in
degree $1$.
\end{lem}
}

\subsubsection{Moduli}

Return to the moduli stack 
$$ \tilde X=\varprojlim_q \tilde X_{\leq q}\,.$$
We have now defined open substacks $\tilde X^s_{\leq q}\subset \tilde
X_{\leq q}$ of stable 
algebras.  We let  $\tilde X^s_q$ the preimage of $\tilde X^s_{\leq q}$ in $\tilde
X$. This is the substack of $q$-stable algebras. Hence we have in 
$\tilde X$ a family of open substacks $\tilde 
X^s_q$, 
parametrized by $q\in\nn$. A point  in $\tilde X$ represents a stable algebra if and
only if it is in almost all open substacks $\tilde X^s_q\subset \tilde
X$. The locus of stable algebras in $\tilde X$ is
$$\tilde X^s=\bigcup_{N\in\nn}\,\bigcap_{q\geq N}\, \tilde X_q^s\,.$$
We see no obvious reason why $\tilde X^s$ should be an open
substack of $\tilde X$.

\begin{rmk}
We have, for every $N\in \nn$ a diagram
\begin{equation}\label{alliso}
\xymatrix{
\bigcap_{q\geq N}\tilde X_q^s\,\ar@{^{(}->}[r]\ar@{^{(}->}[d]  & \tilde
X_N^s\ar@{->>}[r] & \tilde 
X^s_{\leq N}\\
\tilde X^s}\end{equation}
and we find it reasonable, that there should exist dimension vectors
$\vec d$ and integers $N$, for which all arrows in (\ref{alliso}) are
isomorphisms, so that $\tilde X^s=\tilde X^s_{\leq N}$, and $\tilde
X^s$ is a finite type, separated, (in fact quasi-projective)
Deligne-Mumford stack. 

In the commutative case (where $\vec d$ is a numerical polynomial),
the category of graded algebras generated in degree 1, with fixed
Hilbert polynomial, is bounded.  From this it follows that in the
commutative case the corresponding claim $\tilde X^s=\tilde X^s_{\leq
  N}$, for sufficiently large $N$, holds. Lack of suitable persistence
theorems and flattening stratifications currently keep us from
generalizing this result to the non-commutative case. 
\end{rmk}

\begin{defn}
Call a sufficiently ample (meaning that $H^i\big(\cC,\O(n)\big)$
vanishes for $i>0$ and $n>0$), non-commutative projective scheme
$(\cC,\O,s)$ {\bf stable}, if 
$\bigoplus_{n>0} \Gamma(\cC,\O(n))$ is a stable graded algebra.
\end{defn}

\begin{rmk}
It stands to reason,
by analogy with the commutative case, that for certain $\vec d$, the
stack $\tilde X^s$, or an open substack, 
is a moduli stack for stable non-commutative projective
schemes. Further evidence is provided by the deformation theory
arguments from Section~\ref{secone}, which indicate that the derived
deformation theory of a non-commutative projective scheme coincides
with that of its algebra of homogeneous coordinates. 
\end{rmk}


\end{document}